\input amstex
\documentstyle{amsppt}
%
%
\nopagenumbers

\def\cotan{\operatorname{cotan}} 
\def\vp{\operatorname{v.p.}}
\def\negskp{\hskip -2pt}
\pagewidth{360pt}
\pageheight{606pt}
\leftheadtext{R.~A.~Sharipov}
\rightheadtext{On the solutions of weak normality equations \dots}
\topmatter
\title On the solutions of weak normality equations
in multidimensional case.
\endtitle
\author
Ruslan A. Sharipov
\endauthor
\abstract
The system of weak normality equations constitutes a part in the complete
system of norma\-lity equations. Solutions of each of these two systems of
equations are associated with some definite classes of Newtonian
dynamical systems in Riemannian mani\-folds. In this paper for the
case of simplest flat Riemannian manifold $M=\Bbb R^n$ with $n\,{\ssize
\geqslant}\,3$ we show that there exist solutions of weak normality
equations that do not solve complete system of normality equations
in whole. Hence associated classes of Newtonian dynamical systems
do not coincide with each other.
\endabstract
\address Rabochaya street 5, 450003, Ufa, Russia
\endaddress
\email \vtop{\hsize=280pt\noindent
R\_\hskip 1pt Sharipov\@ic.bashedu.ru\newline
ruslan-sharipov\@usa.net\newline
http:/\negskp/www.geocities.com/CapeCanaveral/Lab/5341}
\endemail
\endtopmatter
\loadbold
\document
\head
1. Introduction.
\endhead
     Let $M$ be a Riemannian manifold of the dimension $n$, and let
$S$ be a hypersurface in $M$. One of the ways for deforming $S$ consists
in moving points of $S$ along trajectories of some Newtonian dynamical
system. In local coordinates $x^1,\,\ldots,\,x^n$ in $M$ such system
is given by $n$ ordinary differential equations 
$$
\ddot x^k+\sum^n_{i=1}\sum^n_{j=1}\Gamma^k_{ij}\,\dot x^i\,\dot x^j
=F^k(x^1,\ldots,x^n,\dot x^1,\ldots,\dot x^n),\hskip -2em
\tag1.1
$$
where $k=1,\,\ldots,\,n$. Here $\Gamma^k_{ij}=\Gamma^k_{ij}(x^1,\ldots,
x^n)$ are components of metric connection for basic metric $\bold g$ of
the manifold $M$, and $F^k$ are components of force vector $\bold F$.
They determine force field of dynamical system \thetag{1.1}.\par
     In order to obtain a shift of hypersurface $S$ we set up the
following Cauchy problem for the system of ordinary differential
equations \thetag{1.1}:
$$
\xalignat 2
&\quad x^k\,\hbox{\vrule height 8pt depth 8pt width 0.5pt}_{\,t=0}
=x^k(p),
&&\dot x^k\,\hbox{\vrule height 8pt depth 8pt width 0.5pt}_{\,t=0}=
\nu(p)\cdot n^k(p).\hskip -2em
\tag1.2
\endxalignat
$$
Here $p$ is a point of hypersurface $S$, while $n^k(p)$ are components
of unitary normal vector to $S$. The function $\nu(p)$ is interpreted
as modulus of initial velocity $|\bold v|$; we assume it to be a smooth
function of $p$. For the fixed point $p$ initial data \thetag{1.2}
determine some trajectory of dynamical system \thetag{1.1} coming out
from the point $p$ and being perpendicular to $S$ at this point. Let's
map the point $p$ to the point $p(t)$ on such trajectory. If such
correspondence can be extended to the whole surface $S$, then we have
a shifting map $f_t\!:S\to S_t$. In general case the theorem on smooth
dependence of the solution of ODE upon initial data (see \cite{1} or
\cite{2}) warranties only the possibility to extend this map to some
neighborhood of the point $p$ in $S$. Let $S'$ be such neighborhood.
When $t$ is sufficiently close to zero, shifting maps $f_t\!:S\to S_t$
are diffeomorphisms, \pagebreak their images are smooth hypersurfaces.
In whole,
diffeomorphisms $f_t\!:S'\to S'_t$, which are defined locally, are
glued into a one-parametric family of local\footnotemark\ diffeomorphisms
$f_t\!:S\to S_t$.\footnotetext{Locality here means that domain of the
map $f_t$ depends on $t$. For sufficiently large $t$ it can be empty
at all.}
\definition{Definition 1.1} One-parametric set of local diffeomorphisms
$f_t\!:S\to S_t$ determined by the equations \thetag{1.1} and by initial
data \thetag{1.2} is called {\bf a shift} of $S$ along trajectories of
dynamical system \thetag{1.1}. Such shift is called {\bf a normal shift}
if all hypersurfaces $S_t$ produced by this shift are perpendicular to
its trajectories.
\enddefinition
    Let $p_0$ be some point of hypersurface $S$, and let $\nu_0$ be some
nonzero number. We normalize the function $\nu(p)$ in \thetag{1.2} by the
condition:
$$
\nu(p_0)=\nu_0.\hskip -2em
\tag1.3
$$
\definition{Definition 1.2} Newtonian dynamical system \thetag{1.1}
with force field $\bold F$ is called a system {\bf admitting the normal
shift in strong sense\footnotemark}\ if for any hypersurface $S$ in $M$,
for any point $p_0\in S$, and for any real number $\nu_0\neq 0$ there
exists a neighborhood $S'$ of the point $S$, and there exists smooth
nonzero function $\nu(p)$ in $S'$ normalized by the condition \thetag{1.3}
and such that the shift $f_t\!:S'\to S'_t$ determined by this function
is a normal shift in the sense of definition~1.1.
\enddefinition
\footnotetext{Earlier we used the definition without normalizing
condition \thetag{1.3} for the function $\nu(p)$. Such definition was
called the normality condition. The definition~1.2 strengthens this
condition making it more strict with respect to the choice of force
field $\bold F$ of the dynamical system \thetag{1.1}. It is called
the strong normality condition.}
\adjustfootnotemark{-2}
     Definitions~1.1 and 1.2 underlie in the base of the theory of
dynamical systems admitting the normal shift. This theory was
developed in papers \cite{3--18}; results of these papers were
used in preparing theses \cite{19} and \cite{20}. In paper \cite{8}
we have derived the following equations, which were called
{\bf weak normality equations}:
$$
\cases
\dsize\sum^n_{i=1}\left(v^{-1}\,F_i+\shave{\sum^n_{j=1}}
\tilde\nabla_i\left(N^j\,F_j\right)\right)P^i_k=0,\\
\vspace{2ex}
\aligned
 &\sum^n_{i=1}\sum^n_{j=1}\left(\nabla_iF_j+\nabla_jF_i-2\,v^{-2}
 \,F_i\,F_j\right)N^j\,P^i_k\,+\\
 &+\sum^n_{i=1}\sum^n_{j=1}\left(\frac{F^j\,\tilde\nabla_jF_i}{v}
 -\sum^n_{r=1}\frac{N^r\,N^j\,\tilde\nabla_jF_r}{v}\,F_i\right)P^i_k=0.
 \endaligned
\endcases\hskip -2em
\tag1.4
$$
In paper \cite{9} {\bf additional normality equations} were derived:
$$
\pagebreak
\cases
\aligned
&\sum^n_{i=1}\sum^n_{j=1}P^i_\varepsilon\,P^j_\sigma\left(
\,\shave{\sum^n_{m=1}}N^m\,\frac{F_i\,\tilde\nabla_mF_j}{v}
-\nabla_iF_j\right)=\\
&\quad=\sum^n_{i=1}\sum^n_{j=1}P^i_\varepsilon\,P^j_\sigma\left(
\,\shave{\sum^n_{m=1}}N^m\,\frac{F_j\,\tilde\nabla_mF_i}{v}
-\nabla_jF_i\right),
\endaligned\\
\vspace{2ex}
\dsize\sum^n_{i=1}\sum^n_{j=1}P^j_\sigma\,\tilde\nabla_jF^i\,
P^\varepsilon_i=\sum^n_{i=1}\sum^n_{j=1}\sum^n_{m=1}
\frac{P^j_m\,\tilde\nabla_jF^i\,P^m_i}{n-1}\,P^\varepsilon_\sigma.
\endcases\hskip -2em
\tag1.5
$$
Leaving the equations \thetag{1.4} and \thetag{1.5} with no comments
for a while, now we shall formulate two theorems binding these
equations with definitions~1.1 and 1.2.
\proclaim{Theorem 1.1} Newtonian dynamical system in two-dimensional
Riemannian ma\-nifold $M$ satisfies strong normality condition if and
only if its force field $\bold F$ satisfies weak normality equations
\thetag{1.4} for $v=|\bold v|\neq 0$.
\endproclaim
\proclaim{Theorem 1.2} Newtonian dynamical system in Riemannian manifold
$M$ of the dimension $n\geqslant 3$ satisfies strong normality condition
if and only if its force field $\bold F$ satisfies the normality equations
\thetag{1.4} and \thetag{1.5} for $v=|\bold v|\neq 0$.
\endproclaim
    Theorems~1.1 and 1.2 show that two and multidimensional cases are
substantially different. In multidimensional case complete system of
the equations \thetag{1.4} and \thetag{1.5} is strongly overdetermined,
so that in reducing it we find that it is integrable in explicit form.
General solution for the equations \thetag{1.4} and \thetag{1.5} is
determined by two arbitrary functions $h=h(w)$ and $W=W(x^1,\ldots,x^n,v)$:
$$
F_k(p,\bold v)=\frac{h(W)\,N_k}{W_v}-v\sum^n_{i=1}\frac{\nabla_iW}
{W_v}\,\bigl(2\,N^i\,N_k-\delta^i_k\bigr).\hskip -2em
\tag1.6
$$
Formula \thetag {1.6} was obtained in thesis \cite {19}. Here $\nabla_iW$ 
is the partial derivative of the function $W$ in $i$-th coordinate of the 
point $p$, while $W_v$ is partial derivative of $W$ in the variable $v$, 
which is is interpreted as modulus of velocity vector: $v=|\bold v|$.
\par 
    In two-dimensional case we have only the equations \thetag {1.4}. 
Here they were reduced to one scalar partial differential equation of 
the second order, then various special solutions of this equation 
were constructed (see thesis \cite{20}). Most of these solutions 
correspond to the force fields that can be obtained by means of 
formula \thetag{1.6} when taking $n=2$ in it. These force fields
were called fields of {\bf multidimensional type}. However, one of 
the most important results of thesis \cite{20} is that there
essentially two-dimensional solution of the equations \thetag{1.4}
was constructed. This solution is not given by formula \thetag{1.6}. 
Here we also construct the solution of the equations \thetag{1.4},
which is not expressed by formula \thetag{1.6}, but in multidimensional
case.\par
    Formula \thetag{1.6} obtained in thesis \cite{19} describes all
Newtonian dynamical systems admitting the normal shift in Riemannian
manifolds of the dimension $n\geqslant 3$. However, when thesis
\cite{19} has been already written, new statement of the problem
of normal shift was found. It leads to the equations \thetag{1.4} in
pure form (without additional equations \thetag{1.5}). Indeed, in
papers \cite{3--18} and in theses \cite{19} and \cite{20} we considered
only smooth hypersurfaces $S$ and we choosed sufficiently small
values of $t$ for the sift $f_t\!:S\to S_t$ to result in
smooth hypersurfaces $S_t$ only. If we eliminate this restriction for
$t$, then in the process of shifting we sometimes can observe singular
points on hypersurface $S_t$ (they are called {\bf caustics}). In
particular, under the definite sircumstances hypersurface $S_t$ can
contract into a point at a time for some $t=t_0$. This process is
called {\bf collapse}. Immediately after the callapse for $t>t_0$
we shall observe a blow-up of the point into a series of expanding
hypersurfaces $S_t$. The idea to consider the blow-ups of points by
means of Newtonian dynamical systems were suggested to me by
A.~V.~Bolsinov and A.~T.~Fomenko when I was reporting results of
thesis \cite{19} in the seminar at Moscow State University in
February of 2000. This idea was realized in paper~\cite{21} and in
paper~\cite{22}.\par
    Let $p_0$ be some point of Riemannian manifold $M$. Let's consider
the set of all unitary vectors in tangent space $T_{p_0}(M)$. They can
be interpreted as radius-vectors of the points of unit sphere $\sigma$
in $T_{p_0}(M)$. Let $q\in\sigma$ and let $\bold n(q)$ be the
radius-vector of the point $q$ on unit sphere. Let's fix some constant
number $\nu_0\neq 0$ and set up the following Cauchy problem for the
equations of Newtonian dynamical system \thetag{1.1}:
$$
\xalignat 2
&\quad x^k\,\hbox{\vrule height 8pt depth 8pt width 0.5pt}_{\,t=0}
=x^k(p_0),
&&\dot x^k\,\hbox{\vrule height 8pt depth 8pt width 0.5pt}_{\,t=0}=
\nu_0\cdot n^k(q)\hskip -2em
\tag1.7
\endxalignat
$$\par
\parshape 4 0pt 360pt 0pt 360pt 0pt 360pt 160pt 200pt
\noindent
For the fixed $q$ initial data \thetag{1.2} determine some trajectory of
dynamical system \thetag{1.1} coming out fromn the point $p_0$. To the
point $q\in\sigma$ we put into correspondence the point $p(t)$ on such
trajectory. \vadjust{\vskip 13pt\hbox to 0pt{\kern 5pt\hbox{
}\hss}\vskip -13pt}Theorem on 
smooth dependence of the solution of ODE upon
initial data (see \cite{1} or \cite{2}) says that we can extend this map
to some neighborhood of the point $q$ in $\sigma$. If $\sigma'$ is
such neighborhood of the point $q$, then we have one-parametric family
of diffeomorphisms $f_t\!:\sigma'\to S'_t$. Due to compactness of unit
sphere $\sigma$ we can glue local maps into one map $f_t\!:\sigma\to S_t$.
Here diffeomorphisms $f_t\!:\sigma\to S_t$ are determined globally on the
whole sphere $\sigma$, though parameter $t$ can be restricted by some
interval $(-\varepsilon,\,+\varepsilon)$ on real axis as before.
\definition{Definition 1.3}\parshape 3 160pt 200pt 160pt 200pt 0pt 360pt
One parametric family of diffeomorphisms $f_t\!:\sigma\to S_t$ given by
the equations \thetag{1.1} and initial data \thetag{1.7} is called
{\bf a blow-up} of the point $p_0$  along trajectories of dynamical
system \thetag{1.1}. It is called {\bf a normal blow-up} if all
hypersurfaces $S_t$ arising in this blow-up are perpendicular to
its trajectories.
\enddefinition
\definition{Definition 1.4} Newtonian dynamical system \thetag{1.1}
with force field $\bold F$ in Riemannian manifold $M$ is called {\bf
admitting normal blow-up of points} if for any point $p_0\in M$, and
for any positive constant $\nu_0$ initial data \thetag{1.7} determine
normal blow-up of this point along trajectories of dynamical system
\thetag{1.1}.
\enddefinition
    Definitions~1.3 and 1.4 were first formulated in paper \cite{21}.
They introduced new object: a class of newtonian dynamical systems
admitting the normal blow-up of points in Riemannian manifolds. In
paper \cite{21} it was shown (see theorem~12.1) that this new class
of systems comprises the class of dynamical systems admitting the
normal shift of hypersurfaces, which was previously considered. More
exact description of new class of dynamical systems is given by the
following theorem proved in paper \cite{22}.
\proclaim{Theorem 1.3} Newtonian dynamical system \thetag{1.1} on
Riemannian manifold $M$ \ admits normal blow-up of points if and only
if its force field $\bold F$ satisfies weak normality equations
\thetag{1.4} for $|\bold v|\neq 0$.
\endproclaim
    From theorems~1.1 and 1.3 we see that class of dynamical systems
admitting the normal shift of hypersurfaces in two dimensional case
$n=2$ coincides with the class of systems admitting normal blow-up of
points. This was proved in \cite{22}. But for the multidimensional
case $n\geqslant 3$ the question on {\bf coinciding} or {\bf not
coinciding} of these two classes remained open. I.~A.~Taimanov was
strongly interested in this question during my report in the seminar
of Yu.~G.~Reshetnyak at MI SB RAS (Mathematical Institute of Siberian
\pagebreak
Brunch of Russian Academy of Sciences) in October of 2000. The main goal
of present paper is to give answer to this question and, thus, eliminate
one more obstacle for defending\footnotemark\ thesis \thetag{19}.
\footnotetext{The matter is that in Russia one can pretend for the
degree of Doctor of Sciences only upon writing thesis and passing
so called Defense Procedure in Specialized Council.}
\adjustfootnotemark{-1}
\head
2. Normality equations and extended tensor fields.
\endhead
     Let's consider the force field $\bold F$ in the equations of
Newtonian dynamics \thetag{1.1}. Left hand side of these equations
are components of acceleration vector $\nabla_t\bold v$  (covariant
derivative of velocity vector with respect to parameter $\tau$ along
the trajectory). Therefore $F^k$ are components of tangent vector to
$M$. However, they depend on double set of arguments: on coordinates
$x^1,\,\ldots,\,x^n$ of the point $p\in M$ and on the components
of tangent vector $\bold v\in T_p(M)$. Pair $q=(p,\bold v)$ is a
point of tangent bundle $TM$, so that $p=\pi(q)$. Thus, considering
the equations of the form \thetag{1.1}, we come to the concept of
{\bf extended vector field}. Its generalization is a concept of
{\bf extended tensor field}.
\definition{Definition 2.1} Function $\bold X$ that maps each point
$q=(p,\bold v)$ of tangent bundle $TM$ to some tensor of the type
$(r,s)$ from the space $T^r_s(p,M)$ at the point $p=\pi(q)$ is called
{\bf extended tensor field} of the type $(r,s)$ on the manifold $M$.
\enddefinition
    The concept of extended tensor field stems from thesis \cite{23}
of Finsler, which gave rise to Finslerian geometry. Another approach
to constructing extended tensor fields consists in rising them from
$M$ to $TM$. Here they constitute some special subset in the set of
traditional tensor fields on $TM$. Such tensor fields were considered
in the book \cite{24}, they were called {\bf semibasic tensor fields}.
I am grateful to N.~S.~Dairbekov from IM SB RAS, who noted that
theories of extended and semibasic tensor fields are isomorphic to
each other.\par
     Below we shall use theory of extended tensor fields, which is
described in details in Chapters \uppercase
\expandafter{\romannumeral 2}--\uppercase\expandafter{\romannumeral 4}
of thesis \cite{19}. It is based on the definition~2.1.\par
    For the beginning let's consider some particular examples of extended
tensor fields on the Riemannian manifold $M$.\par
    1. Let's take the point $q=(p,\bold v)$ of tangent bundle $TM$ and
let's map it to the vector $\bold v$ belonging to tangent space $T_p(M)$.
This yields an extended vector field, which is called {\bf the field of
velocity}.\par
    2. Let's map the point $q=(p,\bold v)$ of $TM$ to the number $v=|\bold
v|$. This yields an extended scalar field, which is called {\bf the field
of modulus of velocity vector}.\par
    3. Extended field of unitary vectors $\bold N$ is determined as the
ratio of two previous fields: $\bold N=\bold v/v$.\par
    4. Extended field of operators $\bold P$ is formed by operators of
orthogonal projection onto the hyperplanes perpendicular to velocity
vector $\bold v$. Its components can be written explicitly: $P^i_j=
\delta^i_j-N^i\,N_j$.\par
    Components of all above fields are present in normality equations
\thetag{1.4} and \thetag{1.5}. Moreover, in these equations we see the
operators of covariant differentiation $\nabla$ and $\tilde\nabla$;
they called {\bf spatial} and {\bf velocity} gradients. The simplest
way to define them is to use explicit formulas in coordinates:
$$
\gather
\tilde\nabla_mX^{i_1\ldots\,i_r}_{j_1\ldots\,j_s}=\frac{\partial
X^{i_1\ldots\,i_r}_{j_1\ldots\,j_s}}{\partial v^m},
\hskip -2em
\tag2.1\\
\allowdisplaybreak
\aligned
&\nabla_mX^{i_1\ldots\,i_r}_{j_1\ldots\,j_s}=\frac{\partial
X^{i_1\ldots\,i_r}_{j_1\ldots\,j_s}}{\partial x^m}
-\sum^n_{a=1}\sum^n_{b=1}v^a\,\Gamma^b_{ma}\,\frac{\partial
X^{i_1\ldots\,i_r}_{j_1\ldots\,j_s}}{\partial v^b}\,+\\
&+\sum^r_{k=1}\sum^n_{a_k=1}\!\Gamma^{i_k}_{m\,a_k}\,X^{i_1\ldots\,
a_k\ldots\,i_r}_{j_1\ldots\,\ldots\,\ldots\,j_s}
-\sum^s_{k=1}\sum^n_{b_k=1}\!\Gamma^{b_k}_{m\,j_k}
X^{i_1\ldots\,\ldots\,\ldots\,i_r}_{j_1\ldots\,b_k\ldots\,j_s}.
\endaligned\hskip -2em
\tag2.2
\endgather
$$\medskip
\head
3. Scalar ansatz.
\endhead
    In order to simplify the normality equations \thetag{1.4} and
\thetag{1.5} in paper \cite{18} the scalar ansatz was suggested.
Note that it is applicable either in two-dimensional case $n=2$,
and in multidimensional case $n\geqslant 3$ as well:
$$
F_k=A\,N_k-|\bold v|\,\sum^n_{i=1}P^i_k\,\tilde\nabla_i A.
\hskip -2em
\tag3.1
$$
Scalar ansatz \thetag{3.1} follows from first normality equation
in the system \thetag{1.4}. Substituting \thetag{3.1} back into
this equation, we turn it into identity.\par
    Formula \thetag{3.1} expresses force vector $\bold F$ through
one extended scalar field $A$. This scalar field $A$ can be expressed
back through $\bold F$:
$$
A=\sum^n_{i=1}F^i\,N_i.\hskip -2em
\tag3.2
$$
Substituting \thetag{3.1} into the second equation in the system
\thetag{1.4}, we obtain the following equation for scalar field $A$:
$$
\gathered
\sum^n_{s=1}\left(\nabla_sA+|\bold v|\sum^n_{q=1}\sum^n_{r=1}
P^{qr}\,\tilde\nabla_qA\,\tilde\nabla_r\tilde\nabla_sA\right.-\\
\vspace{0.5ex}
-\left.\shave{\sum^n_{r=1}}N^r\,A\,\tilde\nabla_r\tilde\nabla_sA
-|\bold v|\shave{\sum^n_{r=1}}N^r\,\nabla_r\tilde\nabla_sA\right)
P^s_k=0.
\endgathered\hskip -2em
\tag3.3
$$
Formulas \thetag{3.1} and \thetag{3.2} establishes one-to-one
correspondence between solutions of the equations \thetag{1.4} and
\thetag{3.3}. If force field $\bold F$ satisfies the equations
\thetag{1.4} and \thetag{1.5} simultaneously, then it is expressed
by formula \thetag{1.6}. Scalar field $A$ corresponding to such
force field is given by formula
$$
A=\frac{h(W)}{W_v}-v\sum^n_{i=1}\frac{N^i\,\nabla_iW}{W_v}.\hskip -2em
\tag3.4
$$
Formula \thetag{3.4} admits gauge transformations, which change
functions $h=h(w)$ and $W=W(x^1,\ldots,x^n,v)$, but which don't
change $A$:
$$
\aligned
&W(x^1,\ldots,x^n,v)\longrightarrow \rho(W(x^1,\ldots,x^n,v)),
\vspace{1ex}
&h(w)\longrightarrow h(\rho^{-1}(w))\,\,\rho'(\rho^{-1}(w)).
\endaligned\hskip -2em
\tag3.5
$$
If the function $h(w)$ is nonzero, then by means of gauge transformations
\thetag{3.5} it can be made identically equal to unity (see thesis \cite{19}
and succeeding paper \cite{21}). Therefore, instead of formula \thetag{3.4}
with two arbitrary functions, we can use two formulas with one arbitrary
function:
$$
A=\cases\dsize\frac{1}{W_v}-v\sum^n_{i=1}\frac{N^i\,\nabla_iW}{W_v}
&\text{for \ }h=1,\\
\vspace{1ex}
\dsize-v\sum^n_{i=1}\frac{N^i\,\nabla_iW}{W_v}&\text{for \ }h=0.
\endcases\hskip -2em
\tag3.6
$$
Now the problem stated in section 1 is reformulated as follows: one
should find a solution of the equations \thetag{3.3} that cannot be
expressed by formula \thetag{3.6} neither for $h=1$ nor $h=0$.
\head
4. Spatially homogeneous force field with axial symmetry.
\endhead
    Let $M$ be a space $\Bbb R^n$ with standard Euclidean metric.
We shall construct the required solution of the equations \thetag{3.3}
in this simples case. In the space $\Bbb R^n$ covariant derivatives
$\nabla_m$ and $\tilde\nabla_m$ given by formulas \thetag{2.1} and
\thetag{2.2} turn to partial derivatives: $\nabla_m=\partial/\partial
x^m$ and $\tilde\nabla_m=\partial/\partial v^m$. We restrict our
consideration to spatially homogeneous force fields, in $\Bbb R^n$
they depend only on velocity vector, but don't depend on coordinates
$x^1,\,\ldots,\,x^n$. For the corresponding function $A$ this yields
$A=A(v^1,\ldots,v^n)$. When substituting this function into the
equation \thetag{3.3}, this equation is reduced to the following
one:
$$
\sum^n_{s=1}\left(|\bold v|\sum^n_{q=1}\sum^n_{r=1}
P^{qr}\,\tilde\nabla_qA\,\tilde\nabla_r\tilde\nabla_sA
-\shave{\sum^n_{r=1}}N^r\,A\,\tilde\nabla_r\tilde\nabla_sA
\right)P^s_k=0.\hskip -2em
\tag4.1
$$\par
\parshape 17 0pt 360pt 0pt 360pt 160pt 200pt 160pt 200pt
160pt 200pt 160pt 200pt 160pt 200pt 160pt 200pt 160pt 200pt
160pt 200pt 160pt 200pt 160pt 200pt 160pt 200pt 160pt 200pt
160pt 200pt 160pt 200pt
0pt 360pt\noindent
    Let's mark some direction in$\Bbb R^n$ determined by some constant
unitary vector $\bold m$. Without loss of generality we can assume
\vadjust{\vskip 17pt\hbox to 0pt{\kern 5pt\hbox{
}\hss}\vskip -17pt}$\bold m$ 
to be directed along $n$-th coordinate axis. Let's expand $\bold v$
into a sum 
$$
\bold v=u\cdot\bold n+w\cdot\bold m,\hskip -2em
\tag4.2
$$
where $\bold n\perp\bold m$ and $|\bold n|=1$ All directions perpendicular
to $\bold m$ are assumed to be equivalent. Therefore we choose function $A$
depending only on two variables: $A=A(u,w)$. One should study whether
such choice is compatible with normality equations \thetag{4.1}. For this
purpose we calculate partial derivatives
$$
\tilde\nabla_iA=\frac{\partial A}{\partial v^i}=
\cases A_u\cdot\dfrac{v_i}{u}&\text{for \ }i<n,\\
\vspace{1ex}
A_w&\text{for \ }i=n.\endcases\hskip -2em
\tag4.3
$$
Formula \thetag{4.3} shows that vector of velocity gradient $\tilde
\nabla A$ belong to the linear span of vectors $\bold m$ and $\bold n$
from the expansion \thetag{4.2}:
$$
\tilde\nabla A=A_u\cdot\bold n+A_w\cdot\bold m.\hskip -2em
\tag4.4
$$
Denote by $\bold B$ the vector with the following components:
$$
B^r=\sum^n_{q=1}P^{qr}\,\tilde\nabla_qA.
$$
This is the projection of velocity gradient $\tilde\nabla A$ to the
hyperplane perpendicular to velocity vector. One can obtain explicit
formula for the vector $\bold B$. For this purpose let's denote by
$\theta$ the angle between vectors $\bold m$ and $\bold v$. Then let's
expand projections $\bold P\bold m$ and $\bold P\bold n$ in the base
composed by unitary vectors $\bold m$ and $\bold n$:
$$
\aligned
&\bold P\bold m=\sin^2\theta\cdot\bold m-\sin\theta\,
\cos\theta\cdot\bold n,\\
&\bold P\bold n=\cos^2\theta\cdot\bold n-\sin\theta\,
\cos\theta\cdot\bold m.
\endaligned
\hskip -2em
\tag4.5
$$
From \thetag{4.4} and \thetag{4.5} we derive the following expression
for $\bold B$:
$$
\aligned
\bold B=\,&(A_u\,\cos^2\theta-A_w\,\sin\theta\,\cos\theta)\cdot\bold n\,+\\
&+\,(A_w\,\sin^2\theta-A_u\,\sin\theta\,\cos\theta)\cdot\bold m.
\endaligned
\hskip -2em
\tag4.6
$$
Components of the vector $\bold B$ are present in the equation \thetag{4.1}.
Now we can write this equation in the following form:
$$
\sum^n_{s=1}\sum^n_{r=1}(v\,B^r-A\,N^r)\ \tilde\nabla_r\tilde
\nabla_sA\ P^s_k=0.\hskip -2em
\tag4.7
$$
Now let's calculate the derivatives $\tilde\nabla_r\tilde\nabla_sA$, which
are present in the equation \thetag{4.7}. In order to do it we rewrite
formula \thetag{4.3} as follows:
$$
\tilde\nabla_sA=A_u\,\frac{v_s-m_s\,(\bold m\,|\bold v)}{u}+
A_w\,m_s.
$$
Here $(\bold m\,|\bold v)$ is the scalar product of vectors $\bold m$ and
$\bold v$. Now for second order derivatives $\tilde\nabla_r\tilde\nabla_sA$
by direct calculations we get
$$
\gather
\tilde\nabla_r\tilde\nabla_sA=
A_{uu}\,\frac{v_r-m_r\,(\bold m\,|\bold v)}{u}\cdot
\frac{v_s-m_s\,(\bold m\,|\bold v)}{u}\,+\\
+\,A_{uw}\left(m_r\,\frac{v_s-m_s\,(\bold m\,|\bold v)}{u}+
\frac{v_r-m_r\,(\bold m\,|\bold v)}{u}\,m_s\right)
+A_{ww}\,m_r\,m_s.
\endgather
$$
The above expression is rather complicated. In order to simplify it
let's note that $(\bold m\,|\bold v)=v\,\cos\theta$, $(\bold n\,|
\bold v)=v\,\sin\theta=u$. Hence
$$
\frac{v_r-m_r\,(\bold m\,|\bold v)}{u}=
\frac{n_r\,(\bold n\,|\bold v)}{u}=n_r.
$$
When applied to second order derivatives $\tilde\nabla_r\tilde\nabla_sA$,
this equality yields:
$$
\tilde\nabla_r\tilde\nabla_sA=A_{uu}\,n_r\,n_s
+A_{uw}\,(m_r\,n_s+n_r\,m_s)+A_{ww}\,m_r\,m_s.
$$
Denote $\bold b=\cos\theta\cdot\bold n-\sin\theta\cdot\bold m$. This is
unitary vector belonging to the linear span of vectors $\bold v$ and
$\bold m$ and being perpendicular to $\bold v$. Using this vector we can
rewrite the relationships \thetag{4.5} as follows:
$$
\xalignat 2
&\bold P\bold n=\cos\theta\cdot\bold b,
&&\bold P\bold m=-\sin\theta\cdot\bold b.
\endxalignat
$$
Now we are able to contract $\tilde\nabla_r\tilde\nabla_sA$ with the
components of projector $\bold P$:
$$
\gathered
\sum^n_{s=1}\tilde\nabla_r\tilde\nabla_sA\ P^s_k=
A_{uu}\,\cos\theta\,n_r\,b_k-\,A_{uw}\,\sin\theta\,n_r\,b_k\,+\\
\vspace{1ex}
+\,A_{uw}\,\cos\theta\,m_r\,b_k-A_{ww}\,\sin\theta\,m_r\,b_k.
\endgathered\hskip -2em
\tag4.8
$$
Let's use \thetag{4.8} for to rewrite \thetag{4.7} in more explicit
form:
$$
\gathered
\sum^n_{r=1}(A_{uu}\,\cos\theta-A_{uw}\,\sin\theta)
\,(v\,B^r-A\,N^r)\,n_r\,b_k\,+\\
+\sum^n_{r=1}(A_{uw}\,\cos\theta-A_{ww}\,\sin\theta)
\,(v\,B^r-A\,N^r)\,m_r\,b_k=0.
\endgathered\hskip -3em
\tag4.9
$$
Vector $\bold b$ is nonzero. Its components cannot vanish simultaneously.
Hence in left hand side of \thetag{4.9} we can collect common multiple
$b_k$ and cancel it. Sums in $r$ are scalar products. Therefore we get
$$
\gathered
(A_{uu}\,\cos\theta-A_{uw}\,\sin\theta)
\,(v\,\bold B-A\,\bold N\,|\,\bold n)\,+\\
+\,(A_{uw}\,\cos\theta-A_{ww}\,\sin\theta)
\,(v\,\bold B-A\,\bold N\,|\,\bold m)=0.
\endgathered\hskip -3em
\tag4.10
$$
In order to calculate scalar products in \thetag{4.10} we use the
expansion \thetag{4.6} for the vector $\bold B$ and the expansion
$\bold N=\cos\theta\cdot\bold m+\sin\theta\cdot\bold n$ for the
vector $\bold N$:
$$
\gathered
(A_{uu}\,\cos\theta-A_{uw}\,\sin\theta)
(v\,A_u\,\cos^2\theta\,-\\
-\,v\,A_w\,\sin\theta\,\cos\theta-A\,\sin\theta)
+(A_{uw}\,\cos\theta-A_{ww}\,\sin\theta)\times\\
\times(v\,A_w\,\sin^2\theta-v\,A_u\,\sin\theta\,\cos\theta
-A\,\cos\theta)=0.
\endgathered\hskip -3em
\tag4.11
$$
\proclaim{Theorem 4.1} System weak normality equations \thetag{3.3}
written with respect to the function $A=A(x^1,\ldots,x^n,v^1,\ldots,
v^n)$ in flat Euclidean case $M=\Bbb R^n$ admits the substitution
$A=A(u,w)$, where $u=\sqrt{(v^1)^2+\ldots+(v^{n-1})^2\,}$ and $w=v^n$.
Thereby it is reduced to the single differential equation \thetag{4.11},
where $v=\sqrt{(u)^2+(w)^2\,}$ and $\theta=\arccos(w/v)$.
\endproclaim
    The equation \thetag{4.11} is still rather complicated. In order to
simplify it we transform it to polar coordinates $v$ and $\theta$ in the
plane of variables $u$ and $w$, i\. e\. we do the following change of
variables$$
\xalignat 2
&\quad u=v\,\sin\theta,&&w=v\,\cos\theta.\hskip -3em
\tag4.12
\endxalignat
$$
Let's calculate whether how partial derivatives are transformed under
the change of variables \thetag{4.12}. For the first order derivatives
we have
$$
\pagebreak
\xalignat 2
&A_u=A_v\,\sin\theta+A_\theta\,\frac{\cos\theta}{v},
&&A_w=A_v\,\cos\theta-A_\theta\,\frac{\sin\theta}{v}.
\endxalignat
$$
Then let's calculate second order partial derivatives:
$$
\gather
A_{uu}=A_{vv}\,\sin^2\!\theta+A_{v\theta}\,\frac{\sin{2\theta}}{v}
+A_{\theta\theta}\,\frac{\cos^2\!\theta}{v^2}
+A_v\,\frac{\cos^2\!\theta}{v}-A_\theta\,\frac{\sin2\theta}{v},\\
\vspace{1ex}
A_{uw}=A_{vv}\,\frac{\sin 2\theta}{2}+A_{v\theta}\,\frac{\cos
2\theta}{v}-A_{\theta\theta}\,\frac{\sin 2\theta}{2\,v^2}
-A_v\,\frac{\sin 2\theta}{2\,v}-A_\theta\,\frac{\cos 2\theta}{v^2},\\
\vspace{1ex}
A_{ww}=A_{vv}\,\cos^2\!\theta-A_{v\theta}\,\frac{\sin\theta}{v}+
A_{\theta\theta}\,\frac{\sin^2\!\theta}{v^2}+A_{v}\,
\frac{\sin^2\!\theta}{v}+A_\theta\,\frac{\sin 2\theta}{v^2}.
\endgather
$$
And finally, let's substitute all the above expressions for partial
derivatives into the equation \thetag{4.11}. Thereby the equation
\thetag{4.11} crucially simplifies and takes the form
$$
\frac{A\,A_\theta}{v}+\frac{A_\theta\,A_{\theta\theta}}{v}
+A_\theta\,A_v=A\,A_{v\theta}.\hskip -3em
\tag4.13
$$
\proclaim{Theorem 4.2} System weak normality equations \thetag{3.3}
written with respect to the function $A=A(x^1,\ldots,x^n,v^1,\ldots,
v^n)$ in flat Euclidean case $M=\Bbb R^n$ admits the substitution
$A=A(v,\theta)$, where $v=|\bold v|$ and $\theta=\arccos(v^n/|\bold v|)$.
Thereby it is reduced to the single differential equation \thetag{4.13}.
\endproclaim
Note that the equation \thetag{4.13} do not depend on the dimension of
the space $M=\Bbb R^n$. It holds either in two-dimensional case $n=2$,
and in multidimensional case $n\geqslant 3$ as well. Moreover, this
equation is well known in the theory of dynamical systems admitting the
normal shift (see paper \cite{10} and thesis \cite{20}). In paper
\cite{10} was shown that the equation \thetag{4.13} is integrable in
quadratures. For this purpose it was first transformed to the following
form:
$$
\frac{A_\theta}{A}\cdot\left(\frac{A_\theta}{A}
\right)^{\!\prime}_{\!\theta}-v\cdot\left(\frac{A_\theta}{A}
\right)^{\!\prime}_{\!v}+\frac{A_\theta}{A}+
\left(\frac{A_\theta}{A}\right)^3=0.
$$
This form of the equation \thetag{4.13} says that we should denote
$A_\theta/A=b$. Then for $b=b(v,\theta)$ we obtain quasilinear
partial differential equation of the first order:
$$
b\,b_\theta-v\,b_v+b+b^3=0.\hskip -3em
\tag4.14
$$
In the equation \thetag{4.14} it is convenient to do another one
change of variables, taking $b=\cotan z$, where $z=z(v,\theta)$.
This brings the equation \thetag{4.14} to the form
$$
z_\theta-v\,\frac{\sin z}{\cos z}\,z_v-1=0.\hskip -3em
\tag4.15
$$
The equation \thetag{4.15} with the use of method of characteristics
(see \cite{25}). Characteristics of the equation \thetag{4.15} are
the solutions of the following system of ODEs:
$$
\pagebreak
\cases
\dot\theta=1,\\
\dot z=1,\\
\dot v=-v\,\dfrac{\sin z}{\cos z}.
\endcases\hskip -3em
\tag4.16
$$
Note, that from \thetag{4.16} one can derive the following equalities:
$$
\xalignat 2
&\dot\theta-\dot z=0,
&&\frac{\cos z}{v^2}\,\dot v+\frac{\sin z}{v}\,\dot z=0.
\endxalignat
$$
These equalities can be integrated. They mean that the system of equations
\thetag{4.16} has the pair of first integrals $I_1$ and $I_2$:
$$
\xalignat 2
&I_1=\theta-z,&&I_2=\frac{\cos z}{v}.
\hskip -3em
\tag4.17
\endxalignat
$$
General solution of the equation \thetag{4.15} is determined by first
integrals \thetag{4.17} in implicit form by means of functional
equation
$$
\Phi(I_1,I_2)=0,\hskip -3em
\tag4.18
$$
where $\Phi$ is some arbitrary function of two variables.\par
     Formula \thetag{4.18} proves the integrability of the equation
\thetag{4.13} in quadratures. But in general it yields only local
solution of this equation. Our goal is to construct global solution
of the equation \thetag{4.13}. It should be a smooth function in
direct product of two intervals: closed interval $[0,\,\pi]$ for
the variable $\theta$ and open interval $(0,\,+\infty)$ for the
variable $v$. At the ends of the interval $[0,\,\pi]$ one should
provide the boundary conditions
$$
\xalignat 2
&A_\theta\,\hbox{\vrule height 8pt depth 8pt width
0.5pt}_{\,\theta=0}=0,
&&A_\theta\,\hbox{\vrule height 8pt depth 8pt width
0.5pt}_{\,\theta=\pi}=0.
\hskip -3em
\tag4.19
\endxalignat
$$
They appear because the function $A(v,\theta)$ should correspond to
to the function of several variables $A(v^1,\ldots,v^n)$ with axial
symmetry, being its restriction to the plane passing through the
axis of symmetry.\par
    Now let's proceed with constructing the required solution of the
equation \thetag{4.13}. Note that functional relation of two first
integrals $I_1$ and $I_2$ written as \thetag{4.18} for some cases
can be given by the function of one variable $y=f(w)$. Let's write
\thetag{4.18} as $I_1=f(I_2)$, i\. e\. let's consider the following
functional equation:
$$
\theta-z=f\!\left(\frac{\cos z}{v}\right).\hskip -3em
\tag4.20
$$
For $f(w)$ we choose smooth increasing function with decreasing
derivative; we assume that $f(w)$ is defined in semiopen interval
$[0,\,+\infty)$ and $f(0)=\pi/2$. Let's also assume that $f(w)$
is not restricted and grows to infinity as $w\to+\infty$. Its
graph is shown on Fig\.~4.2. By means of function $y=f(w)$ we
construct a family of functions $y=F_{\sssize[\ssize v\sssize]}(z)$
depending on $v$ as parameter:
$$
y=F_{\sssize[\ssize v\sssize]}(z)=f\!\left(\frac{\cos z}{v}
\right).\hskip -3em
\tag4.21
$$
Graphs of the functions are shown on Fig\.~4.3. We use them in order
to solve the equation \thetag{4.20} graphically. For this purpose we
consider a family of straight lines being graphs of the following
functions:
$$
y=F_{\sssize[\theta]}(z)=\theta-z.\hskip -3em
\tag4.22
$$
Denote by $v_{\sssize\min}$ the value of the derivative $f'(w)$ at
the point $w=0$:
$$
v_{\sssize\min}=f'(w)\,\hbox{\vrule height 8pt depth 8pt width
0.5pt}_{\,w=0}=0.
$$
Suppose that the values of parameters $v$ and $\theta$ satisfy the
following inequalities:
$$
\xalignat 2
&v_{\sssize\min}<v<+\infty,
&&0\leqslant\theta\leqslant\pi.
\hskip -3em
\tag4.23
\endxalignat
$$
From Fig.~4.2 and Fig.~4.3 we see that for these values of parameters
graphs of functions \thetag{4.21} and \thetag{4.22} intersect at
unique point and determine smooth function $z=z(v,\theta)$ satisfying
the equation \thetag{4.15}. \vadjust{\vskip 12pt\hbox to 0pt{\kern 5pt
\hbox{
}\hss}\vskip 160pt}For each fixed value of $v$ from the domain determined
by inequalities \thetag{4.23} the function $z(v,\theta)$ is increasing
in $\theta$, it takes all values from closed interval $[-\pi/2,\,+\pi/2]$:
$$
\xalignat 2
&z\,\hbox{\vrule height 8pt depth 8pt width
0.5pt}_{\,\theta=0}=-\frac{\pi}{2},
&&z\,\hbox{\vrule height 8pt depth 8pt width
0.5pt}_{\,\theta=\pi}=\frac{\pi}{2}.
\hskip -3em
\tag4.24
\endxalignat
$$
Exactly at one point $\theta_0=\theta_0(v)$ in the interval $[0,\,\pi]$
this function vanishes, while its derivative in $\theta$ at this point
is equal to unity:
$$
z_\theta(v,\theta_0(v))=1.\hskip -3em
\tag4.25
$$
Note also that any fixed value of $\theta$ the function $z(v,\theta)$
is increasing function in $v$, though the interval of its values here
is more narrow, and it depends on $\theta$.\par
    Function $z(v,\theta)$, which is constructed graphically, determines
the function $b=\cotan z=b(v,\theta)$. For the fixed value of $v$ it is
decreasing function in $\theta$, but it has an infinite break at the
point $\theta_0=\theta_0(v)$. From the equality \thetag{4.25} we derive
$$
b(v,\theta)=\frac{1}{\theta-\theta_0}+O(1)\text{\ \ for \ }\theta
\to\theta_0.\hskip -3em
\tag4.26
$$
From \thetag{4.24} we obtain that $b(v,\theta)$ vanishes at both ends
of interval $[-\pi/2,\,+\pi/2]$:
$$
\xalignat 2
&b\,\hbox{\vrule height 8pt depth 8pt width
0.5pt}_{\,\theta=0}=0,
&&b\,\hbox{\vrule height 8pt depth 8pt width
0.5pt}_{\,\theta=\pi}=0.
\hskip -3em
\tag4.27
\endxalignat
$$
Using $b(v,\theta)$, now we define the function $A(v,\theta)$ by
the following formula
$$
A(v,\theta)=\exp\left(\vp\!\!\!\int\limits^{\,\,\theta}_{-\pi/2}
\!\!b(v,\tau)\,d\tau\right)\!.\hskip -3em
\tag4.28
$$
From \thetag{4.26} it follows that the function \thetag{4.28} vanishes
at the point $\theta_0=\theta_0(v)$, while from \thetag{4.27} we derive
boundary conditions \thetag{4.19} for this function. By construction
the function \thetag{4.28} is a solution of the equation \thetag{4.13}.
However its domain \thetag{4.23} do not embrace the whole phase space.
In order to expand its domain one should note that if $A(v,\theta)$ is
the solution of the equation \thetag{4.13}, then the product $C(v)\cdot
A(v,\theta)$, where $C=C(v)$ is an arbitrary smooth function, is also
the solution of this equation. Let's choose the function $C(v)$ such
that
$$
\aligned
&C(v)=0\text{\ \ for \ }v\leqslant v_{\sssize\min},\\
&C(v)=1\text{\ \ for \ }v\geqslant v_0>v_{\sssize\min}.
\endaligned
$$
Now, instead of formula \thetag{4.28}, we define the function
$A(v,\theta)$ by formula 
$$
A(v,\theta)=C(v)\cdot\exp\left(\vp\!\!\!\int\limits^{\,\,\theta}_{-\pi/2}
\!\!b(v,\tau)\,d\tau\right)\!.\hskip -3em
\tag4.29
$$
Function \thetag{4.29} is determined in all phase space, except for
those points, where $v=|\bold v|=0$. It is the solution of the equation
\thetag{4.13} and it satisfies boundary conditions \thetag{4.19}.
\head
5. Theorem on non-coincidence of classes.
\endhead
    The solution of the equation \thetag{4.13} constructed by formula
\thetag{4.29} determines the function $A(v^1,\ldots,v^n)$, which, in
turn, determines force field $\bold F$ of Newtonian dynamical system
in $\Bbb R^n$ admitting the normal blow-up of points. The construction
of this field has functional arbitrariness due to the function $f=f(w)$
in \thetag{4.20} and the function $C=C(v)$ in \thetag{4.29}. Formula
\thetag{3.6} for scalar field $A$ corresponding to dynamical systems
admitting the normal shift of hypersurfaces also has functional
arbitrariness due to the function $W=W(x^1,\ldots,x^n,v)$. Let's study
which part of this arbitrariness remains if we assume $A$ to be
spatially homogeneous function with axially symmetric dependence on
$\bold v$. The dependence on the direction of velocity vector $\bold v$
in \thetag{3.6} is completely determined by the sum
$$
B=v\sum^n_{i=1}\frac{N^i\,\nabla_iW}{W_v}=\sum^n_{i=1}
\frac{v^i\,\nabla_iW}{W_v}=\frac{(\bold v\,|\,\nabla W)}{W_v}
\hskip -2em
\tag5.1
$$
Changing $v$ by $-v$, we change the sign of this sum, but the value
of $v=|\bold v|$, which is the argument of function $W$, remains
unchanged. For the case $h=1$ this yields
$$
\frac{1}{W_v}=\frac{A(\bold v)+A(-\bold v)}{2}.\hskip -2em
\tag5.2
$$
In spatially homogeneous case right hand side of \thetag{5.2} doesn't
depend on coordinates $x^1,\,\ldots,\,x^n$. Hence for $h=1$ the
quantity $1/W_v$ depends only on $v$. Let's denote it by $H(v)$. For
the quantity $B$ this yields
$$
B=\cases A(\bold v)-H(v)&\text{for \ }h=1,\\
\ A(\bold v)&\text{for \ }h=0.\endcases
\hskip -2em
\tag5.3
$$
Thus, for both cases the quantity $B=B(\bold v)$ in \thetag{5.3} do not
depend on coordinates $x^1,\ldots,x^n$. Let $\bold e_1,\,\ldots,\,\bold
e_n$ be unitary vectors directed along coordinate axes. Substituting
$\bold v=v\cdot\bold e_i$ into the sum \thetag{5.1}, we get 
$$
m_i=\frac{\nabla_iW}{W_v}=\frac{B(v\cdot\bold e_i)}{v}.\hskip -2em
\tag5.4
$$
The quantities $m_i$ in \thetag{5.4} determine some vector $\bold m$.
Due to \thetag{5.3} they do not depend on $x^1,\,\ldots,\,x^n$. Hence
$\bold m=\bold m(v)$. Substituting the quantities \thetag{5.4} into
the sum \thetag{5.1} and further into the formula \thetag{3.6} for
scalar field $A$, we get
$$
A=\cases H(v)+(\bold v\,|\,\bold m)&\text{for \ }h=1,\\
\quad(\bold v\,|\,\bold m)&\text{for \ }h=0.\endcases\hskip -2em
\tag5.5
$$
{\it {\bf Conclusion}: the condition of spatial homogeneity reduces
functional arbitrariness in \thetag{3.6} to the choice of $(n+1)$
functions of one variable. These are the function $H(v)$ and components
of the vector $\bold m(v)$ in formula \thetag{5.5}}.\par
     For the fixed value of $v=|\bold v|$ the function $A(\bold v)$ in
\thetag{5.5} possess axial symmetry with the axis directed along the
vector $\bold m(v)$. While for the function \thetag{4.29} the axis of
symmetry doesn't depend on $v$, it is directed along the vector
$\bold e_n$. Aiming to express the function \thetag{4.29} by formula
\thetag{5.5}, we should choose
$$
\bold m(v)=\frac{C(v)}{v}\cdot\bold e_n.
$$
In variables $v$ and $\theta$ in two-dimensional plane of axial section
this yields
$$
A=\cases H(v)+C(v)\,\cos\theta &\text{for \ }h=1,\\
\quad C(v)\,\cos\theta &\text{for\ }h=0.\endcases\hskip -2em
\tag5.6
$$
Similar to \thetag{4.29}, formula \thetag{5.6} contain functional
arbitrariness determined by two functions of one variable. But this
functional arbitrariness does not affect the dependence of $A$ upon
angular variable $\theta$. While the dependence on $\theta$ in
formula \thetag{4.29} is much more complicated. It is determined by
the choice of function $f(w)$ in the equation \thetag{4.20}; in
general case it is not reduced to trigonometric function $y=\cos
\theta$. Thus, formula \thetag{4.29} determines some solution of the
system of weak normality equations \thetag{1.4}, which is not the
solution for additional normality equations \thetag{1.5}. Therefore
we can formulate the main result of present paper.
\proclaim{Theorem 5.1} In multidimensional case $n\geqslant 3$ class
of Newtonian dynamical \nolinebreak systems admitting normal blow-up
of points is the expansion of the class of systems \nolinebreak
admitting the normal shift of hypersurfaces, and it doesn't coincide
\pagebreak with the latter one.
\endproclaim
     Note that in thesis \cite{20} for two-dimensional case $n=2$
Andrey Boldin has constructed another (more explicit) solution
for the equation \thetag{4.13}. It is expressed through
elliptic functions. However, it is local and it doesn't satisfy the
conditions \thetag{4.19}. Therefore this solution cannot be used in
multidimensional case $n\geqslant 3$ for proving theorem~5.1.
\head
6. Acknowledgements.
\endhead
    I am grateful to Yu.~G.~Reshetnyak and I.~A.~Taimanov for the
invitation to visit Novosibirsk and for the opportunity to report results
of thesis \cite{19} and succeeding papers \cite{21} and \cite{22} in two
seminars at Mathematical Institute of SB RAS (Siberian Brunch of Russian
Academy of Sciences). I am grateful to I.~A.~Taimanov, N.~S.~Dairbekov,
and B.~A.~Sharafutdinov for special attention to the subject of my
reports.\par
    I am also grateful to S.~P.~Tsarev and O.~V.~Kaptsov for the
invitation to visit Krasnoyarsk and for opportunity to make report
in the seminar at the Institute of Computational Modeling of SB RAS.
This allowed me to continue my journey along Trans-Siberian main
line and reach most eastern point on the Earth I have ever been before.
The most western point is Chicago, where I have been in 1996 thanks
to invitation from S.~I.~Pinchuk and E.~Bedford (Indiana University).
I am grateful to S.~P.~Tsarev for useful remark concerning my thesis
\cite{19} and for the interest to thesis \cite{20} of my student
A.~Yu\.~Boldin.\par
    I am grateful to all participants of seminars mentioned above for
fruitful discussions, which stimulated preparing this paper.\par
     This work is supported by grant from Russian Fund for Basic Research
(project No\nolinebreak\.~00\nolinebreak-01-00068, coordinator 
Ya\.~T.~Sultanaev), and by grant from Academy of Sciences of the
Republic Bashkortostan (coordinator N.~M.~Asadullin). I am grateful
to these organizations for financial support.

\enddocument
\end